\documentclass[11pt,a4paper,reqno]{amsart}
\usepackage{graphicx}
\usepackage{subfigure}
 \theoremstyle{definition}
 
 \theoremstyle{remark}
 
\begin{document}
\setcounter{page}{1}
\begin{flushleft}
\end{flushleft}
\bigskip
\bigskip
\title[G. Zhou, Y. Zhou, B. Ahmad, a. Alsaedi: Finite difference/element method ... ] {Finite difference/element method for time-fractional Navier-Stokes equations}
\author[Appl. Comput. Math., V.xx, N.xx,  20xx]{Guang-an Zou$^{1}$,  Yong Zhou$^{2,3}$, Bashir Ahmad$^{3}$, Ahmed Alsaedi$^{3}$}
\thanks{$^1$School of Mathematics and Statistics, Henan University, Kaifeng 475004, P. R. China\\
$^2$Faculty of Mathematics and Computational Science, Xiangtan University\\Xiangtan,  Hunan 411105, P.R. China\\
$^3$Nonlinear Analysis and Applied Mathematics (NAAM) Research Group, Faculty of Science\\
King Abdulaziz University, Jeddah 21589, Saudi Arabia\\
\indent\,\,\,e-mail: zouguangan@henu.edu.cn, yzhou@xtu.edu.cn, bashirahmad\_qau@yahoo.com, aalsaedi@hotmail.com
\\ \indent
  \em \,\,\,Manuscript received xx}

\begin{abstract}
We apply a composite idea of semi-discrete finite difference approximation in time and  Galerkin finite element method in space to solve the Navier-Stokes equations with Caputo derivative of order $0<\alpha<1$. The stability properties and convergence error estimates for both the semi-discrete and fully discrete schemes are obtained. Numerical example is provided to illustrate the validity of theoretical results.


\end{abstract}
\maketitle

\smallskip
\section{Introduction}In this paper, we study the following Navier-Stokes equations with time-fractional derivative in a bounded subset of $\Omega\subset \mathbf{R}^{2}$ with a smooth boundary $\partial \Omega$:
\begin{align*}
\begin{cases}
^{C}D_{t}^{\alpha}u+(u\cdot \nabla)u-\nu\Delta u+\nabla p=f,~~~\mathrm{div}~u=0,~\forall~(x,t)\in \Omega\times(0,T],\\
~u(x,t)|_{\partial \Omega}=0,~t\in (0,T],\\
~u(x,0)=u_{0},~x\in \Omega,\\
\end{cases} \tag{1.1}
\end{align*}
where $^{C}D_{t}^{\alpha}$ represents the Caputo-type fractional derivative of order $\alpha\in(0,1)$, $u=(u_{1}(x,t),$ $ u_{2}(x,t))$ denotes the velocity field at a point $x\in\Omega$ and $t\in[0,T]$, $\nu>0$ is viscosity coefficient, $p=p(x,t)$ represents the pressure field, $f=f(x,t)$ is the external force and $u_{0}=u_{0}(x)$ is the initial velocity.

Notice that the problem (1.1) reduces to the classical Navier-Stokes equations (NSEs) for $\alpha=1$. The existence and non-existence of solutions for the NSEs have been discussed in [1].
Chemin et al.[2] studied the global regularity for the large solutions to the NSEs. Miura [3] focused on the uniqueness of mild solutions to the NSEs. Germain [4] presented the uniqueness criteria for the solutions of the Cauchy problem associated to the NSEs. The existence of global weak solutions for supercritical NSEs was discussed [5]. The lower bounds on blow up solutions for the NSEs in homogeneous Sobolev spaces were studied in [6]. The numerical methods for solving the NSEs have been investigated by many authors [7,8,9,10,11,12,23]. The study of time-fractional Navier-Stokes equations (TFNSEs) has become a hot topic of research due to its significant role in simulating the anomalous diffusion in fractal media. There are also some analytical methods available for solving the TFNSEs. Momani and Odibat [13] applied Adomian decomposition method to obtain the analytical solution of the TFNSEs. In [14, 15], the homotopy perturbation (transform) method was used to find the analytical solution of the TFNSEs. Wang and Liu [16] solved TFNSEs by applying the transform methods. Concerning the existence of global and local mild solutions to TFNSEs, see Carvalho-Neto and Gabriela [17], Zhou and Peng [18]. Moreover, Zhou and Peng [19] investigated the existence of weak solutions and optimal control for TFNSEs, while Peng et al.[20] presented the rigorous exposition of local solutions of TFNSEs in Sobolev space. However, one can notice that there are only a few works related to the numerical solution of the TFNSEs. The details of meshless local Petrov-Galerkin method based on moving Kriging interpolation for solving the TFNSEs can be found in the literature [21]. The purpose of this paper is to present finite difference/element method to obtain the numerical solution of TFNSEs.

The rest of the paper is arranged as follows. In Section 2, we give some notations and preliminaries. Section 3 deals with a semi-discrete scheme for the TFNSEs, which is based on a mixed finite element method in space. We also discuss the stability and error estimates of this semi-discrete scheme. In Section 4, we use a finite difference approximation to discrete time direction to obtain the fully discrete scheme. The stability and error estimates for the discrete schemes are also found. In Section 5, numerical results are discussed to confirm our theoretical analysis. Conclusions are given in the final section.
\bigskip
\section{Notations and preliminaries}

In this section, we present some preliminary concepts of the functional spaces. Firstly, we introduce the following Hilbert spaces:
\begin{align*}
X=H_{0}^{1}(\Omega)^{2},~Y=L^{2}(\Omega)^{2},~M=L_{0}^{2}(\Omega)=\{v\in L^{2}(\Omega);\int_{\Omega}v dx=0\},
\end{align*}
where the space $L^{2}(\Omega)$ is associated with the usual inner product $(\cdot,\cdot)$ and the norm $\|\cdot\|$. The space $X$ is associated with the following inner product and equivalent norm:
\begin{align*}
((u,v))=(\nabla u,\nabla v),~\|u\|_{X}=\|\nabla u\|_{0}=\| u\|_{1}.
\end{align*}

Denote by $V$ and $H$ the closed subsets of $X$ and $Y$ respectively, which are  given by
\begin{align*}
V=\{v\in X;\mathrm{div}~ v=0\},~H=\{v\in Y;\mathrm{div}~ v=0,v\cdot n|_{\partial\Omega}=0\}.
\end{align*}

We denote the Stokes operator by $A=-P\triangle$, in which $P$ is the $L^{2}$-orthogonal projection of $Y$ onto $H$. The domain of $A$ is $\mathcal{D}(A)=H^{2}(\Omega)^{2}\cap V$ and let $H^{s}=\mathcal{D}(A^{\frac{s}{2}})$ with the norm $\|v\|_{s}=\|A^{\frac{s}{2}} v\|$. Observe that $H^{2}=\mathcal{D}(A)$, $H^{1}=V$ and $H^{0}=H$.

Next, we define the Riemann-Liouville fractional integral operator of order $\beta$ ($\beta\geq0$) as (see [22])
\begin{align*}
I_{t}^{\beta}g(t)=\frac{1}{\Gamma(\beta)}\int_{0}^{t}(t-s)^{\beta-1}g(s)ds,t>0, \tag{2.1}
\end{align*}
with $I_{t}^{0}g(t)=g(t)$.

The Caputo-type derivative of order $\alpha \in (0,1]$, $^{C}D_{t}^{\alpha}$ in (1.1) is defined by
\begin{align*}
^{C}D_{t}^{\alpha}g(t)=\frac{d}{dt}\{I_{t}^{1-\alpha}[g(t)-g(0)]\}=\frac{d}{dt}\{\frac{1}{\Gamma(1-\alpha)}\int_{0}^{t}(t-s)^{-\alpha}[g(t)-g(0)]ds\}.\tag{2.2}
\end{align*}
Further, the operator $^{C}D_{t}^{-\alpha}$ is defined as
\begin{align*}
^{C}D_{t}^{-\alpha}g(t)=I_{t}^{\alpha}g(t)=\frac{1}{\Gamma(\alpha)}\int_{0}^{t}(t-s)^{\alpha-1}g(s)ds,t>0.\tag{2.3}
\end{align*}
where $\Gamma(\cdot)$ stands for the gamma function $\Gamma(x)=\int_{0}^{\infty}t^{x-1}e^{-t}dt$.

Next we introduce the following continuous bilinear forms $a(\cdot,\cdot)$ and $d(\cdot,\cdot)$ on $X\times X$ and $X\times M$ respectively as follows:
\begin{align*}
a(u,v)=\nu(\nabla u,\nabla v),~u,v\in X, ~d(v,q)=(q,\mathrm{div}v),~v\in X, q \in M,
\end{align*}
and the trilinear form $b(u,v,w)$ on $X\times X\times X$ is given by
\begin{align*}
b(u,v,w)=((u\cdot \nabla)v+\frac{1}{2}(\mathrm{div}u)v,w)=\frac{1}{2}((u\cdot \nabla)v,w)-\frac{1}{2}((u\cdot \nabla)w,v),~u,v,w\in X.
\end{align*}

It is well-known that the trilinear form $b(u,v,w)$ has the following properties:
\begin{align*}
&b(u,v,w)=-b(u,w,v),~b(u,v,v)=0,~u,v,w\in X,\\
&|b(u,v,w)|\leq \mu_{0}\|u\|_{1}\|v\|_{1}\|w\|_{1},~u,v,w\in X.
\end{align*}

In terms of the above notations, the weak formulation of problem (1.1) is as follows: find $(u,p)\in(X,M)$ for all $t\in [0,T]$ such that for all $(v,q)\in(X,M)$:
\begin{align*}
\begin{cases}
^{C}D_{t}^{\alpha}(u,v)+a(u,v)+b(u,u,v)-d(v,p)+d(u,q)=(f,v),\\
u(0)=u_{0}.
\end{cases} \tag{2.4}
\end{align*}

In [19], Zhou and Peng discussed the existence and uniqueness of weak solutions for the problem (2.4). The objective of the present work is to obtain the numerical solution of the problem at hand.
\section{Finite element method for space discretization}

Let $T^{h}(\Omega)=\{K\}$ be a mesh of $\Omega$ with a mesh size function $h(x)$, which is the diameter $h_{K}$ of element $K$ containing $x$.
Assuming $h=h_{\Omega}=\max\limits_{x\in\Omega}h(x)$ be the largest mesh size of $T^{h}(\Omega)$, we introduce the mixed finite element subspace $(X_{h},M_{h})$ of $(X,M)$ and define the subspace $V_{h}$ of $X_{h}$ as
\begin{align*}
V_{h}=\{v_{h}\in X_{h};d(v_{h},q_{h})=0,~\forall q_{h}\in M_{h}\}.
\end{align*}

Let $P_{h}:Y\rightarrow V_{h}$ denote the $L^{2}$-orthogonal projection defined by
\begin{align*}
(P_{h}v,v_{h})=(v,v_{h}),~\forall v\in Y,~v_{h}\in V_{h}\}.
\end{align*}

With the above notations, we need some further basic assumptions on the mixed finite element spaces (Refs.[8,10,12]).

(A1) \textit{Approximation}. For each $(v,q)\in (\mathcal{D}(A),M\cap H^{1}(\Omega))$, there exist approximations $(\pi_{h}v,\rho_{h}q)\in (X_{h},M_{h})$ such that
\begin{align*}
\|v-\pi_{h}v\|_{1}\leq Ch\|A v\|,~\|q-\rho_{h}q\|\leq Ch\|q\|_{1}.\tag{3.1}
\end{align*}

(A2) \textit{Inverse~estimate}. For any $(v,q)\in (X_{h},M_{h})$, the following relations hold:
\begin{align*}
\|\nabla v\|\leq Ch^{-1}\| v\|,~\|q\|\leq Ch^{-1}\|q\|_{-1}.\tag{3.2}
\end{align*}

(A3) \textit{Stability~property}. For any $(v,q)\in (X_{h},M_{h})$, the well-known inf-sup condition holds:
\begin{align*}
\sup\limits_{v\in X_{h}}\frac{d(v,q)}{\|v\|_{1}}\geq \lambda \|q\|,\tag{3.3}
\end{align*}
where $\lambda>0$ is a constant.

Further, the following classical properties hold:
\begin{align*}
&\|v-P_{h}v\|+h\|\nabla (v-P_{h}v)\|\leq Ch^{2}\|A v\|,v\in\mathcal{D}(A);\tag{3.4}\\
&\|v-P_{h}v\|\leq Ch\|\nabla (v-P_{h}v)\|,v\in X.\tag{3.5}
\end{align*}

The standard finite element Galerkin approximation for (2.4) holds as follows: Find $(u_{h},p_{h})\in(X_{h},M_{h})$ for all $t\in [0,T]$ such that for all $(v_{h},q_{h})\in(X_{h},M_{h})$, we have
\begin{align*}
\begin{cases}
^{C}D_{t}^{\alpha}(u_{h},v_{h})+a(u_{h},v_{h})+b(u_{h},u_{h},v_{h})-d(v_{h},p_{h})+d(u_{h},q_{h})=(f,v_{h}),\\
u_{h}(0)=u_{0h}=P_{h}u_{0}.
\end{cases} \tag{3.6}
\end{align*}

With the above semi-discrete approximation, a discrete analogue of the Stokes operator $A$ is defined as $A_{h}=-P_{h}\triangle_{h}$ via the condition $(-\triangle_{h}u_{h},v_{h})=((u_{h},v_{h}))$ for all $u_{h},v_{h}\in X_{h}$. The trilinear form $b(u_{h},v_{h},w_{h})$ satisfies the following properties:
\begin{align*}
&b(u_{h},v_{h},w_{h})=-b(u_{h},w_{h},v_{h}),~b(u_{h},v_{h},v_{h})=0,~u_{h},v_{h},w_{h}\in X_{h};\tag{3.7}\\
&|b(u_{h},v_{h},w_{h})|\leq \mu_{0}\|u_{h}\|_{1}\|v_{h}\|_{1}\|w_{h}\|_{1},~u_{h},v_{h},w_{h}\in X_{h}.\tag{3.8}
\end{align*}

\textbf{Theorem 3.1.} For any $t\in [0,T]$ and $0<\alpha < 1$, let $u_{h}$ be the solution of equation (3.6). Then there exists a positive constant $C$ such that
\begin{align*}
\|u_{h}\|^{2}+\nu_{1}\int_{0}^{t}\|u_{h}\|_{1}^{2}ds\leq \|u_{0h}\|^{2}+\frac{(1-\alpha_{1})T^{1+\beta}}{2\nu(1+\beta)\Gamma(\alpha)}+\frac{\alpha_{1}T}{2\nu\Gamma(\alpha)}\max\limits_{t\in [0,T]}\|f\|_{-1}^{\frac{2}{\alpha_{1}}},
\end{align*}
where $\nu_{1}=\frac{\nu T^{\alpha-1}}{2\Gamma(\alpha)}>0$ is constant.

\textbf{Proof.} Taking $v_{h}=u_{h}$, $q_{h}=p_{h}$ in (3.6) and  using the Young's inequality, we get
\begin{align*}
^{C}D_{t}^{\alpha}\|u_{h}\|^{2}+\nu\|u_{h}\|_{1}^{2}=(f,u_{h}) \leq \|f\|_{-1}\|u_{h}\|_{1}
 \leq \frac{1}{2\nu}\|f\|_{-1}^{2}+\frac{\nu}{2}\|u_{h}\|_{1}^{2},\tag{3.9}
\end{align*}
where $\|f\|_{-1}=\|A^{-\frac{1}{2}}f\|$ denotes the dual operator in $\mathcal{D}(A^{-\frac{1}{2}})$.

Applying the integral operator (2.3) to both sides of (3.9) and using the Young's inequality, we obtain
\begin{align*}
&\|u_{h}\|^{2}+\frac{\nu}{2\Gamma(\alpha)}\int_{0}^{t}(t-s)^{\alpha-1}\|u_{h}\|_{1}^{2}ds\leq\|u_{0h}\|^{2}+\frac{1}{2\nu\Gamma(\alpha)}\int_{0}^{t}(t-s)^{\alpha-1}\|f\|_{-1}^{2}ds\\
&\leq \|u_{0h}\|^{2}+\frac{1-\alpha_{1}}{2\nu\Gamma(\alpha)}\int_{0}^{t}(t-s)^{\frac{\alpha-1}{1-\alpha_{1}}}ds+\frac{\alpha_{1}}{2\nu\Gamma(\alpha)}\int_{0}^{t}\|f\|_{-1}^{\frac{2}{\alpha_{1}}}ds\\
&\leq \|u_{0h}\|^{2}+\frac{(1-\alpha_{1})T^{1+\beta}}{2\nu(1+\beta)\Gamma(\alpha)}+\frac{\alpha_{1}T}{2\nu\Gamma(\alpha)}\max\limits_{t\in [0,T]}\|f\|_{-1}^{\frac{2}{\alpha_{1}}},
\end{align*}
where $\beta=\frac{\alpha-1}{1-\alpha_{1}}$ with $0<\alpha_{1} < 1$.

In view of the inequality
\begin{align*}
\frac{\nu}{2\Gamma(\alpha)}\int_{0}^{t}(t-s)^{\alpha-1}\|u_{h}\|_{1}^{2}ds\geq\frac{\nu T^{\alpha-1}}{2\Gamma(\alpha)}\int_{0}^{t}\|u_{h}\|_{1}^{2}ds,
\end{align*}
it follows that
\begin{align*}
\|u_{h}\|^{2}+\frac{\nu T^{\alpha-1}}{2\Gamma(\alpha)}\int_{0}^{t}\|u_{h}\|_{1}^{2}ds\leq \|u_{0h}\|^{2}+\frac{(1-\alpha_{1})T^{1+\beta}}{2\nu(1+\beta)\Gamma(\alpha)}+\frac{\alpha_{1}T}{2\nu\Gamma(\alpha)}\max\limits_{t\in [0,T]}\|f\|_{-1}^{\frac{2}{\alpha_{1}}}.
\end{align*}

This completes the proof.

\textbf{Theorem 3.2.} For any $t\in [0,T]$, $0<\alpha < 1$, let $(u,p)$ and $(u_{h},p_{h})$ be the solutions of equations (2.4) and (3.6) respectively. Then there exists a positive constant $C$ such that
\begin{align*}
\|u-u_{h}\|\leq Ch^{2},~\|p-p_{h}\|\leq Ch. \tag{3.10}
\end{align*}

\textbf{Proof.} Setting $(\xi,\eta)=(u-u_{h},p-p_{h})$, we deduce from (2.4) and (3.6) that
\begin{align*}
\begin{cases}
^{C}D_{t}^{\alpha}(\xi,v)+a(\xi,v)+b(\xi,u_{h},v)+b(u_{h},\xi,v)+b(\xi,\xi,v)-d(v,\eta)+d(\xi,q)=0,\\
\xi_{0}=u_{0}-P_{h}u_{0}.
\end{cases} \tag{3.11}
\end{align*}

Taking $v=\xi$ and $q=\eta$ in (3.11), we get
\begin{align*}
^{C}D_{t}^{\alpha}\|\xi\|^{2}+\nu\|\xi\|_{1}^{2}+b(\xi,u_{h},\xi)=0.
\end{align*}

Using the properties of $b(u_{h},v_{h},w_{h})$ together with Young's inequality, we obtain
\begin{align*}
^{C}D_{t}^{\alpha}\|\xi\|^{2}+\nu\|\xi\|_{1}^{2}&=b(\xi,\xi,u_{h})\\
&\leq C_{0}\|\xi\|\|\xi\|_{1}\|u_{h}\|_{1}\\
&\leq \nu\|\xi\|_{1}^{2}+C_{1}\|\xi\|^{2}\|u_{h}\|_{1}^{2}.\tag{3.12}
\end{align*}

Applying (2.3) to both sides of (3.12), we have
\begin{align*}
\|\xi\|^{2}\leq \|\xi_{0}\|^{2}+\frac{C_{1}}{\Gamma(\alpha)}\int_{0}^{t}(t-s)^{\alpha-1}\|\xi\|^{2}\|u_{h}\|_{1}^{2}ds.\tag{3.13}
\end{align*}

By means of the generalized integral version of Gronwall's lemma [24], we get
\begin{align*}
\|\xi\|^{2}&\leq \|\xi_{0}\|^{2}\exp(\frac{C_{1}}{\Gamma(\alpha)}\int_{0}^{t}(t-s)^{\alpha-1}\|u_{h}\|_{1}^{2}ds)\\
&\leq \|u_{0}-P_{h}u_{0}\|^{2}\exp[C_{2}(\|u_{0h}\|^{2}+\frac{(1-\alpha_{1})T^{1+\beta}}{2\nu(1+\beta)\Gamma(\alpha)}+\frac{\alpha_{1}T}{2\nu\Gamma(\alpha)}\max\limits_{t\in [0,T]}\|f\|_{-1}^{\frac{2}{\alpha_{1}}})]\\
&\leq Ch^{4}.\tag{3.14}
\end{align*}

Furthermore, setting $v=\xi$ and $q=0$ in (3.11), and using inf-sup condition (3.3), combining (3.1)-(3.5),(3.8),(3.14) and using the integral operator (2.3) in (3.15), we conclude that
\begin{align*}
\|\eta\| &\leq \sup\limits_{V_{h}}\frac{|d(\xi,\eta)|}{\lambda\|\xi\|_{1}}\\
&=\sup\limits_{V_{h}}\frac{|^{C}D_{t}^{\alpha}\|\xi\|^{2}+\nu\|\xi\|_{1}^{2}+b(\xi,u_{h},\xi)|}{\lambda\|\xi\|_{1}}\\
&\leq Ch. \tag{3.15}
\end{align*}

This completes the proof.

\section{Finite difference method for time discretization}

The discretization of time-fractional derivative can be found in [25-32] and references therein.
Here, we will introduce a uniform grid by discretizing the temporal domain $[0,T]$ given by the points: $t_{n}=n\tau$ for $n=0,1,\ldots,N$, with the time-step size $\tau=T/N$. Hence, the Riemann-Liouville fractional integral operator of order $\alpha$ can be discretized as follows:
\begin{align*}
I_{t}^{\alpha}g(t_{n})&=\frac{1}{\Gamma(\alpha)}\sum\limits_{k=1}^{n}\int_{t_{k-1}}^{t_{k}}(t_{n}-s)^{\alpha-1}g(s)ds\\
&=\frac{1}{\Gamma(\alpha)}\sum\limits_{k=1}^{n}(\int_{t_{k-1}}^{t_{k}}(t_{n}-s)^{\alpha-1}g(t_{k})ds)+\gamma_{\alpha}^{n}\\
&=\frac{\tau^{\alpha}}{\Gamma(\alpha+1)}\sum\limits_{k=1}^{n}g(t_{k})[(n-k+1)^{\alpha}-(n-k)^{\alpha}]+\gamma_{\alpha}^{n}\\
&=\frac{\tau^{\alpha}}{\Gamma(\alpha+1)}\sum\limits_{k=0}^{n-1}w_{k}^{\alpha}g(t_{n-k})+\gamma_{\alpha}^{n},\tag{4.1}
\end{align*}
where $w_{k}^{\alpha}=(k+1)^{\alpha}-k^{\alpha}$ and the truncation error $\gamma_{\alpha}^{n}$ is given by
\begin{align*}
\gamma_{\alpha}^{n}&=\frac{1}{\Gamma(\alpha)}\sum\limits_{k=1}^{n}\int_{t_{k-1}}^{t_{k}}(t_{n}-s)^{\alpha-1}[g(s)-g(t_{k})]ds\\
&=\frac{1}{\Gamma(\alpha)}\sum\limits_{k=1}^{n}\int_{t_{k-1}}^{t_{k}}(t_{n}-s)^{\alpha-1}g^{\prime}(\zeta)(s-t_{k})ds,~s<\zeta<t_{k}.
\end{align*}

Therefore, we have
\begin{align*}
|\gamma_{\alpha}^{n}|&\leq\frac{\tau}{\Gamma(\alpha)}\max\limits_{0\leq t\leq t_{k}}|g^{\prime}(t)|\sum\limits_{k=1}^{n}\int_{t_{k-1}}^{t_{k}}(t_{n}-s)^{\alpha-1}ds\\
&\leq \frac{N^{\alpha}\tau^{\alpha+1}}{\Gamma(\alpha+1)}\max\limits_{0\leq t\leq t_{k}}|g^{\prime}(t)|.
\end{align*}

\textbf{Lemma 4.1.} (see \cite{25}) If $g(t)\in C^{1}[0,T]$, then
\begin{align*}
I_{t}^{\alpha}g(t_{n})=\frac{\tau^{\alpha}}{\Gamma(\alpha+1)}\sum\limits_{k=0}^{n-1}w_{k}^{\alpha}g(t_{n-k})+\gamma_{\alpha}^{n},  \tag{4.2}
\end{align*}
where $|\gamma_{\alpha}^{n}|\leq C\tau^{\alpha+1},~n=0,1,\ldots,N$.

\textbf{Lemma 4.2.} (see \cite{25}) For $0< t_{n}\leq t_{N}=T$ and $\alpha>0$, let the coefficient $w_{k}^{\alpha}$ be given by (4.1). Then
\begin{align*}
&(\mathrm{i})~ w_{0}^{\alpha}=1,w_{k}^{\alpha}>0,k=0,1,2,\cdots;\\
&(\mathrm{ii})~ w_{k}^{\alpha}>w_{k+1}^{\alpha},k=0,1,2,\cdots;\\
&(\mathrm{iii})\sum\limits_{k=0}^{n-1}w_{k}^{\alpha}=n^{\alpha}\leq N^{\alpha}.
\end{align*}

Applying the integral operator (2.3) to both sides of (3.6), we obtain
\begin{align*}
&(u_{h},v_{h})+\frac{1}{\Gamma(\alpha)}\int_{0}^{t}(t-s)^{\alpha-1}[a(u_{h},v_{h})+b(u_{h},u_{h},v_{h})-d(v_{h},p_{h})+d(u_{h},q_{h})]ds\\
&=(u_{0h},v_{h})+\frac{1}{\Gamma(\alpha)}\int_{0}^{t}(t-s)^{\alpha-1}(f,v_{h})ds.  \tag{4.3}
\end{align*}

Let $u_{h}^{n}$ and $p_{h}^{n}$ be the numerical solutions of $u_{h}(t)$ and $p_{h}(t)$ at $t=t_{n}$ respectively. By (4.2) and (4.3), our full discrete scheme of equation (2.4) can be defined by seeking
$(u_{h}^{n},p_{h}^{n})\in(X_{h},M_{h})$ such that for all $(v_{h},q_{h})\in(X_{h},M_{h})$:
\begin{align*}
&(u_{h}^{n},v_{h})+\beta_{0}\sum\limits_{k=0}^{n-1}w_{k}^{\alpha}[a(u_{h}^{n-k},v_{h})+b(u_{h}^{n-k},u_{h}^{n-k},v_{h})-d(v_{h},p_{h}^{n-k})+d(u_{h}^{n-k},q_{h})]\\
&=(u_{h}^{0},v_{h})+\beta_{0}\sum\limits_{k=0}^{n-1}w_{k}^{\alpha}(f^{n-k},v_{h}), \tag{4.4}
\end{align*}
where $\beta_{0}=\frac{\tau^{\alpha}}{\Gamma(\alpha+1)}$.

\textbf{Theorem 4.1.} For any $0<\tau<T$, the full discrete scheme (4.4) is unconditionally stable, and that
\begin{align*}
\|u_{h}^{n}\|^{2}+\beta_{1}\|u_{h}^{n}\|_{1}^{2}\leq C^{\dag}(\|u_{h}^{0}\|^{2}+\sum\limits_{k=0}^{n}\|f^{k}\|_{-1}^{2}).
\end{align*}

\textbf{Proof.} Setting $n=1$ in (4.4), we get
\begin{align*}
(u_{h}^{1},v_{h})+\beta_{0}[a(u_{h}^{1},v_{h})+b(u_{h}^{1},u_{h}^{1},v_{h})-d(v_{h},p_{h}^{1})+d(u_{h}^{1},q_{h})]=(u_{h}^{0},v_{h})+\beta_{0}(f^{1},v_{h}). \tag{4.5}
\end{align*}

Taking $v_{h}=u_{h}^{1}$ and $q_{h}=p_{h}^{1}$ in (4.5), we have
\begin{align*}
(u_{h}^{1},u_{h}^{1})+\beta_{0}a(u_{h}^{1},u_{h}^{1})=(u_{h}^{0},u_{h}^{1})+\beta_{0}(f^{1},u_{h}^{1}).
\end{align*}

Making use of Young's inequality, we obtain
\begin{align*}
\|u_{h}^{1}\|^{2}+\beta_{0}\nu\|u_{h}^{1}\|_{1}^{2}\leq \frac{1}{2}[\|u_{h}^{0}\|^{2}+\|u_{h}^{1}\|^{2}]+[\frac{\beta_{0}}{2\nu}\|f^{1}\|_{-1}^{2}+\frac{\beta_{0}\nu}{2}\|u_{h}^{1}\|_{1}^{2}],
\end{align*}
that is,
\begin{align*}
\|u_{h}^{1}\|^{2}+\beta_{1}\|u_{h}^{1}\|_{1}^{2}\leq \frac{1}{2}\|u_{h}^{0}\|^{2}+\frac{\beta_{0}}{2\nu}\|f^{1}\|_{-1}^{2}.
\end{align*}

Assuming $v_{h}=u_{h}^{j}\in X_{h}$ and $q_{h}=p_{h}^{j}\in M_{h}$, the following inequality holds
\begin{align*}
\|u_{h}^{j}\|^{2}+\beta_{1}\|u_{h}^{j}\|_{1}^{2}\leq C(\|u_{h}^{0}\|^{2}+\sum\limits_{k=1}^{j}\|f^{k}\|_{-1}^{2}),~j=2,3,\ldots,n-1.
\end{align*}

Setting $v_{h}=u_{h}^{n-k}$ and $q_{h}=p_{h}^{n-k}$ in (4.4), we get
\begin{align*}
(u_{h}^{n},u_{h}^{n-k})+\beta_{0}\sum\limits_{k=0}^{n-1}w_{k}^{\alpha}a(u_{h}^{n-k},u_{h}^{n-k})=(u_{h}^{0},u_{h}^{n-k})+\beta_{0}\sum\limits_{k=0}^{n-1}w_{k}^{\alpha}(f^{n-k},u_{h}^{n-k}).
\end{align*}

By the elementary identity $ab=\frac{1}{2}(a^{2}+b^{2})-\frac{1}{2}(a-b)^{2}$ and the Young's inequality, we have
\begin{align*}
&\frac{1}{2}[\|u_{h}^{n}\|^{2}+\|u_{h}^{n-k}\|^{2}]+\beta_{0}\nu\sum\limits_{k=0}^{n-1}w_{k}^{\alpha}\|u_{h}^{n-k}\|_{1}^{2}\\
&=(u_{h}^{0},u_{h}^{n-k})+\frac{1}{2}\|u_{h}^{n}-u_{h}^{n-k}\|^{2}
+\beta_{0}\sum\limits_{k=0}^{n-1}w_{k}^{\alpha}(f^{n-k},u_{h}^{n-k})\\
&\leq \frac{1}{2}[\|u_{h}^{0}\|^{2}+\|u_{h}^{n-k}\|^{2}]+\frac{1}{2}\|u_{h}^{n}-u_{h}^{n-k}\|^{2}+\sum\limits_{k=0}^{n-1}w_{k}^{\alpha}(\frac{\beta_{0}}{2\nu}\|f^{n-k}\|_{-1}^{2}+\frac{\beta_{0}\nu}{2}\|u_{h}^{n-k}\|_{1}^{2}).
\end{align*}

Together with Lemma 4.2 (ii) (that is, $0<w_{k+1}^{\alpha}<w_{k}^{\alpha}<1$) and $\frac{\beta_{0}\nu}{2}\sum\limits_{k=1}^{n-1}w_{k}^{\alpha}\|u_{h}^{n-k}\|_{1}^{2}\geq0$, we obtain
\begin{align*}
\|u_{h}^{n}\|^{2}+\beta_{1}\|u_{h}^{n}\|_{1}^{2} &\leq \|u_{h}^{0}\|^{2}+\|u_{h}^{n}-u_{h}^{n-k}\|^{2}
+\frac{\beta_{0}}{2\nu}\sum\limits_{k=0}^{n-1}w_{k}^{\alpha}\|f^{n-k}\|_{-1}^{2}\\
&\leq C^{\dag}(\|u_{h}^{0}\|+\sum\limits_{k=0}^{n}\|f^{k}\|_{-1}^{2}). \tag{4.6}
\end{align*}

The proof is completed.

\textbf{Lemma 4.3.} Let $\nu>0$ be the viscosity coefficient and that
\begin{align*}
\|u_{h}^{0}\|^{2}+\sum\limits_{k=0}^{n}\|f^{k}\|_{-1}^{2}\leq \frac{\beta_{1}\nu}{C^{\dag}\mu_{0}}. \tag{4.7}
\end{align*}
Then
\begin{align*}
a(v_{h},v_{h})+b(v_{h},u_{h}^{n},v_{h})\geq0, \tag{4.8}
\end{align*}
where $\mu_{0}>0$ is defined by (3.8) and $C^{\dag}>0$ is constant.

\textbf{Proof.} Making use of (4.6) and (4.7), we get
\begin{align*}
\|u_{h}^{n}\|_{1}^{2}\leq \frac{\nu}{\mu_{0}},
\end{align*}
that is,
\begin{align*}
\nu-\mu_{0}\|u_{h}^{n}\|_{1}^{2}\geq0.
\end{align*}

By the property of $b(v_{h},u_{h}^{n},v_{h})$, we get
\begin{align*}
a(v_{h},v_{h})+b(v_{h},u_{h}^{n},v_{h})\geq (\nu-\mu_{0}\|u_{h}^{n}\|_{1}^{2})\|v_{h}\|_{1}^{2}\geq 0.
\end{align*}

The proof of the lemma is completed.

\textbf{Theorem 4.2.} For $0<\alpha < 1$, let $(u_{h}(t_{n}),p_{h}(t_{n}))$ and $(u_{h}^{n},p_{h}^{n})$ be the solutions of equations (3.6) and (4.4) respectively. There exists a constant $C$ such that
\begin{align*}
\|u_{h}(t_{n})-u_{h}^{n}\|\leq C \tau^{\alpha+1},~\|p_{h}(t_{n})-p_{h}^{n}\|\leq C \tau^{\alpha+1}. \tag{4.8}
\end{align*}

\textbf{Proof.} Let $\xi^{n}=u_{h}(t_{n})-u_{h}^{n}$ and $\eta^{n}=p_{h}(t_{n})-p_{h}^{n}$. Using (4.2)-(4.4) and noting $\xi^{0}=0$, we deduce
\begin{align*}
&(\xi^{n},v_{h})+\beta_{0}\sum\limits_{k=0}^{n-1}w_{k}^{\alpha}[a(\xi^{n-k},v_{h})+b(\xi^{n-k},u_{h}^{n-k},v_{h})-d(v_{h},\eta^{n-k})+d(\xi^{n-k},q_{h})]\\
&=(\gamma_{\alpha}^{n},v_{h}). \tag{4.9}
\end{align*}

For $n=1$, taking $v_{h}=\xi^{1}$ and $q_{h}=\eta^{1}$ in (4.9), we have
\begin{align*}
(\xi^{1},\xi^{1})+\beta_{0}[a(\xi^{1},\xi^{1})+b(\xi^{1},u_{h}^{1},\xi^{1})]=(\gamma_{\alpha}^{1},\xi^{1}).
\end{align*}

By Cauchy-Schwarz inequality and Lemma 4.3, we get
\begin{align*}
\|\xi^{1}\|\leq \|\gamma_{\alpha}^{1}\|\leq C \tau^{\alpha+1}.  \tag{4.10}
\end{align*}

Let us assume that $\|\xi^{m}\|\leq C \tau^{\alpha+1}$ for $m=2,3,\ldots,n-1$. In order to show that the first inequality in (4.8) holds for $m=n$, we set $v_{h}=\xi^{n-k}$ and $q_{h}=\eta^{n-k}$ in (4.9). Then
\begin{align*}
(\xi^{n},\xi^{n-k})+\beta_{0}\sum\limits_{k=0}^{n-1}w_{k}^{\alpha}[a(\xi^{n-k},\xi^{n-k})+b(\xi^{n-k},u_{h}^{n-k},\xi^{n-k})]=(\gamma_{\alpha}^{n},\xi^{n-k}).
\end{align*}

In view of the elementary identity $ab=\frac{1}{2}(a^{2}+b^{2})-\frac{1}{2}(a-b)^{2}$, Young's inequality and Lemma 4.3, we get
\begin{align*}
\frac{1}{2}[\|\xi^{n}\|^{2}+\|\xi^{n-k}\|^{2}]\leq \frac{1}{2}\|\xi^{n}-\xi^{n-k}\|^{2}+\frac{1}{2}[\|\xi^{n-k}\|^{2}+\|\gamma_{\alpha}^{n}\|^{2}], \tag{4.11}
\end{align*}
which implies that
\begin{align*}
\|\xi^{n}\|^{2}\leq C \tau^{2\alpha+2}. \tag{4.12}
\end{align*}

By inverse estimate (3.2) together  with (4.12), we have
\begin{align*}
\|\xi^{n}\|_{1}\leq Ch^{-1}\|\xi^{n}\| \leq C \tau^{\alpha+1}. \tag{4.13}
\end{align*}

On the other hand, setting $v_{h}=\xi^{1}$ and $q_{h}=0$ for $n=1$ in (4.9) and  making use of Cauchy-Schwarz inequality together with (3.1)-(3.4), (3.8) and (4.13), we get
\begin{align*}
\|\eta^{1}\|&\leq \sup\limits_{V_{h}}\frac{|d(\xi^{1},\eta^{1})|}{\lambda\|\xi^{1}\|_{1}}=\sup\limits_{V_{h}}\frac{|\|\xi^{1}\|^{2}+\beta_{0}[\nu\|\xi^{1}\|_{1}^{2}+b(\xi^{1},u_{h}^{1},\xi^{1})]-(\gamma_{\alpha}^{1},\xi^{1})|}{\lambda\|\xi^{1}\|_{1}}\\
&\leq C \tau^{\alpha+1}. \tag{4.14}
\end{align*}

Using the assumption  $\|\eta^{m}\|\leq C \tau^{\alpha+1}$ for $m=2,3,\ldots,n-1$ and taking $v_{h}=\xi^{n-k}$ and $q_{h}=0$ in (4.9), by (4.12) and (4.13), similar to the derivation of (4.14) for $m=n$, we obtain
\begin{align*}
\|\eta^{n}\|&\leq \sup\limits_{V_{h}}\frac{|d(\xi^{n},\eta^{n})|}{\lambda\|\xi^{n}\|_{1}}\leq \sup\limits_{V_{h}}\frac{|\sum\limits_{k=0}^{n-1}w_{k}^{\alpha}d(\xi^{n-k},\eta^{n-k})|}{\lambda\|\xi^{n}\|_{1}}\\
&=\sup\limits_{V_{h}}\frac{|(\xi^{n},\xi^{n-k})+\beta_{0}\sum\limits_{k=0}^{n-1}w_{k}^{\alpha}[\nu\|\xi^{n-k}\|_{1}^{2}+b(\xi^{n-k},u_{h}^{n-k},\xi^{n-k})]-(\gamma_{\alpha}^{n},\xi^{n-k})|}{\lambda\|\xi^{n}\|_{1}}\\
&\leq C \tau^{\alpha+1}.
\end{align*}

This completes the proof.

Next we give the error estimate for fully discrete scheme.

\textbf{Theorem 4.3.} For $0<\alpha < 1$, let $(u(t_{n}),p(t_{n}))$ and $(u_{h}^{n},p_{h}^{n})$ be the solutions of equations (2.4) and (4.4) respectively. Then there exists a positive constant $C$ such that
\begin{align*}
\|u(t_{n})-u_{h}^{n}\|\leq C(h^{2}+\tau^{\alpha+1}),~\|p(t_{n})-p_{h}^{n}\|\leq C(h+\tau^{\alpha+1}). \tag{4.15}
\end{align*}

\textbf{Proof.} It is easy to show that (4.15) follows from Theorem 3.2 and Theorem 4.2 via triangle inequality.

\section{Numerical example}

In this section, we demonstrate the effectiveness of our numerical methods with the aid of examples. We use mixed finite element method for the discretization of spatial direction and finite difference approximation for time discretization. The convergence rates of numerical solutions with respect to space step $h$ and time step $\tau$ are discussed. We consider the regular (uniform) domain $\Omega=(0,1)\times(0,1)$ and the time interval is chosen to be $[0,1]$ with  the viscosity coefficient $\nu=1.5$.
For an appropriate body force $f$, the analytical solution $(u,p)=((u_{1},u_{2}),p)$ of the unstable flow problem with homogeneous boundary conditions becomes
\begin{align*}
&u_{1}=2x^{2}(x-1)^{2}y(y-1)(2y-1)e^{-t},~u_{2}=-2y^{2}(y-1)^{2}x(x-1)(2x-1)e^{-t},\\
&p=(x^{2}-y^{2})e^{-t},
\end{align*}
which automatically satisfy the initial and boundary conditions.

The errors $\|e^{n}\|$ are computed in $L^{2}$-discrete norm. The results of numerical experiments are compared with analytical solution by the rates of the convergence, which are approximately by
\begin{align*}
\mathrm{Rate }= |\frac{\mathrm{ln}(\|e_{f}^{n}\|/\|e_{c}^{n}\|)}{\mathrm{ln}(N_{f}/N_{c})}|,
\end{align*}
where $\|e_{f}^{n}\|$ and $\|e_{c}^{n}\|$ denote the error on finer grid and coarser grid, $N_{f}$ and  $N_{c}$ represent the numbers of meshes on finer grid and coarser grid, respectively.

The spatial convergence rates for the components of velocity $(u_{1},u_{2})$ and pressure $p$ with fixed time step $\tau=1/8$ with different values of $\alpha$ are shown in Fig.1. The convergence rates of velocity $(u_{1},u_{2})$ are in accordance with spatial convergence order $\mathcal{O}(h^{2})$ and the pressure $p$ are closer to order $\mathcal{O}(h)$. Fig.2 give the temporal convergence rates for the components of velocity and pressure with fixed spatial step $h=1/15$ with different values of $\alpha$. We can see that the rates of convergence are closer to the theoretical convergence order $\mathcal{O}(\tau^{\alpha+1})$.

Fig.3 depicts the numerical solutions of the components of velocity $(u_{1},u_{2})$ and pressure $p$, with $h=1/15$ and $\tau=0.1$, when $\alpha=0.4$ and $\alpha=0.8$, respectively. It is not difficult to find that a pair of warm- and cold-core eddies emerge in the velocity field.

\bigskip
\section{Conclusion}

In this study, the finite difference/element method is presented to solve the TFNSEs and the convergence error estimates for the discrete schemes in $L^{2}$-norm are obtained. We present the numerical experiment to illustrate the accuracy of schemes, and the result fully verify the convergence theory. The numerical examples also confirm the thesis [18,19,20] that in
procedure of citing and novelty of the obtained results. Furthermore, the presented methods and analytical techniques in this work can also be extended to other nonlinear time-fractional partial differential equations.

\section{Acknowledgment}
Guang-an Zou is supported by National Nature Science Foundation of China (Grant No. 11626085), Yong Zhou is supported by National Nature Science Foundation of China (Grant No. 11671339).

\section*{References}

[1]  Lemari\'{e}-Rieusset, P.G. Recent developments in the Navier-Stokes problem,  Chapman  Hall/CRC Research Notes in Mathematics, 431. Chapman  Hall/CRC, Boca Raton, FL, 2002, 395 p.

[2]  Chemin, J.Y.,   Gallagher, I., Paicu, M.  Global regularity for some classes of large solutions to the Navier-Stokes equations, Ann. of Math. (2), V.173, N.2, 2011, pp.983-1012.

[3] Miura, H. Remark on uniqueness of mild solutions to the Navier-Stokes equations, J. Funct. Anal., V.218, N.1, 2005, pp.110-129.

[4] Germain, P. Multipliers, paramultipliers, and weak-strong uniqueness for the Navier-Stokes equations, J. Differential Equations, V.226, N.2, 2006, pp.373-428.

[5] Nahmod, A.R., Pavlovic N., Staffilani, G. Almost sure existence of global weak solutions for supercritical Navier-Stokes equations, SIAM J. Math. Anal. V.45, N.6, 2013, pp.3431-3452.

[6] Robinson, J.C., Sadowski, W.,  Silva, R.P.  Lower bounds on blow up solutions of the three-dimensional Navier-Stokes equations in homogeneous Sobolev spaces, J. Math. Phys., V.53, N.11,  2012, 115618, 15 pp.

[7] Ingram, R. A new linearly extrapolated Crank-Nicolson time-stepping scheme for the Navier-Stokes equations, Math. Comp., V.82, N.284, 2013, pp.1953-1973.

[8] Bernardi, C.,  Raugel, G.A.  conforming finite element method for the time-dependent Navier-Stokes equations, SIAM J. Numer. Anal., V.22, N.3, 1985, pp.455-473.

[9] He, Y., Sun, W. Stability and convergence of the Crank-Nicolson/Adams-Bashforth scheme for the time-dependent Navier-Stokes equations, SIAM J. Numer. Anal., V.45, N.2, 2007, pp.837-869.

[10] Girault, V., Raviart, P.A. Finite element methods for Navier-Stokes equations. Theory and algorithms, Springer-Verlag, Berlin, 1986, 374 p.

[11] Kaya, S.,  Rivi\`{e}re, B.  A discontinuous subgrid eddy viscosity method for the time-dependent Navier-Stokes equations, SIAM J. Numer. Anal., V.43, N.4, 2005, pp.1572-1595.

[12] Shan, L.,  Hou, Y. A fully discrete stabilized finite element method for the time-dependent Navier-Stokes equations, Appl. Math. Comput., V.215, N.1, 2009, pp.85-99.

[13] Momani, S., Odibat, Z. Analytical solution of a time-fractional Navier-Stokes equation by Adomian decomposition method, Appl. Math. Comput., V.177, N.2, 2006, pp.488-494.

[14] Ganji, Z.Z., Ganji, D.D., Ganji, Ammar D., Rostamian, M.  Analytical solution of time-fractional Navier-Stokes equation in polar coordinate by homotopy perturbation method, Numer. Methods Partial Differential Equations, V.26, N.1, 2010, pp.117-124.

[15] Kumar, D.,  Singh, J.,  Kumar, S. A fractional model of Navier-Stokes equation arising in unsteady flow of a viscous fluid, J. Ass. Arab Univ. Basic Appl. Sci., V.17, 2015, pp.14-19.

[16] Wang, K., Liu, S. Analytical study of time-fractional Navier-Stokes equation by using transform methods, Adv. Differential Equ., N.61, 2016, pp.12.

[17]  De Carvalho-Neto, P.M., Gabriela, P.  Mild solutions to the time fractional Navier-Stokes equations in $R^{N}$, J. Differential Equations, V.259, N.7, pp.2948-2980.

[18] Zhou, Y., Peng, L. On the time-fractional Navier-Stokes equations,  Comput. Math. Appl., V.73, N.6, 2017, pp.874-891.

[19] Zhou, Y., Peng, L. Weak solutions of the time-fractional Navier-Stokes equations and optimal control, Comput. Math. Appl., V.73, N.6, 2017, pp.1016-1027.

[20] Peng, L.,  Zhou, Y., Ahmad, B.,   Alsaedi, A. The Cauchy problem for fractional Navier-Stokes equations in Sobolev spaces, Chaos Solitons Fractals, V.102, 2017, pp.218-228.

[21] Thamareerat, N., Luadsong, A., Aschariyaphotha, N.  The meshless local Petrov-Galerkin method based on moving Kriging interpolation for solving the time fractional Navier-Stokes equations, SpringerPlus,  V.5, N.417, 2016, pp.19.

[22] Kilbas, A.A.,  Srivastava, H.M., Trujillo, J.J. Theory and applications of fractional differential equations, Elsevier, 2006, 523 p.

[23]  Layton, W.J., Labovschii, A.,   Manica, C.C., Neda, M., Rebholz, L. G. The stabilized, extrapolated
trapezoidal finite element method for the Navier-Stokes equations,  Comput. Methods Appl. Mech. Eng., V.198,  2009, pp.958-974.

[24] Kruse R. Strong and weak approximation of semilinear stochastic evolution equations,  Springer, Cham, 2014, 177 p.

[25] Zeng, F., Li, C., Liu, F., Turner, I. The use of finite difference/element approaches for solving the time-fractional subdiffusion equation, SIAM J. Sci. Comput., V.35, N.6, 2013, pp.A2976-A3000.

[26] Stynes, M.,  O'Riordan, E.,  Gracia, J.L.  Error analysis of a finite difference method on graded meshes for a time-fractional diffusion equation, SIAM J. Numer. Anal., V.55, N.2, pp.2017, 1057-1079.

[27] Stynes, M.,  Gracia, J.L. Preprocessing schemes for fractional-derivative problems to improve their convergence rates, Appl. Math. Lett., V.74, 2017, pp.187-192.

[28] Kopteva, N.,  Stynes, M. Analysis and numerical solution of a Riemann-Liouville fractional derivative two-point boundary value problem, Adv. Comput. Math.,  V.43, N.1, 2017, 77-99.

[29] Zeng, F., Li, C.,  Liu, F.,  Turner, I. Numerical algorithms for time-fractional subdiffusion equation with second-order accuracy, SIAM J. Sci. Comput., V.37, N.1, 2015, pp.A55-A78.

[30] Zheng, M., Liu, F.,  Liu, Q.,  Burrage, K.,   Simpson, M.J. Numerical solution of the time fractional reaction-diffusion equation with a moving boundary, J. Comput. Phys.,  V.338, 2017, pp.493-510.

[31] Cui, M.R. Compact alternating direction implicit method for two-dimensional time fractional diffusion equation, J. Comput. Phys., V.231, N.6, 2012, pp.2621-2633.

[32] Jiang, Y.J.,  Ma, J.T.  High-order finite element methods for time-fractional partial differential
equations, J. Comput. Appl. Math., V.235, N.11, 2011, 3285-3290.




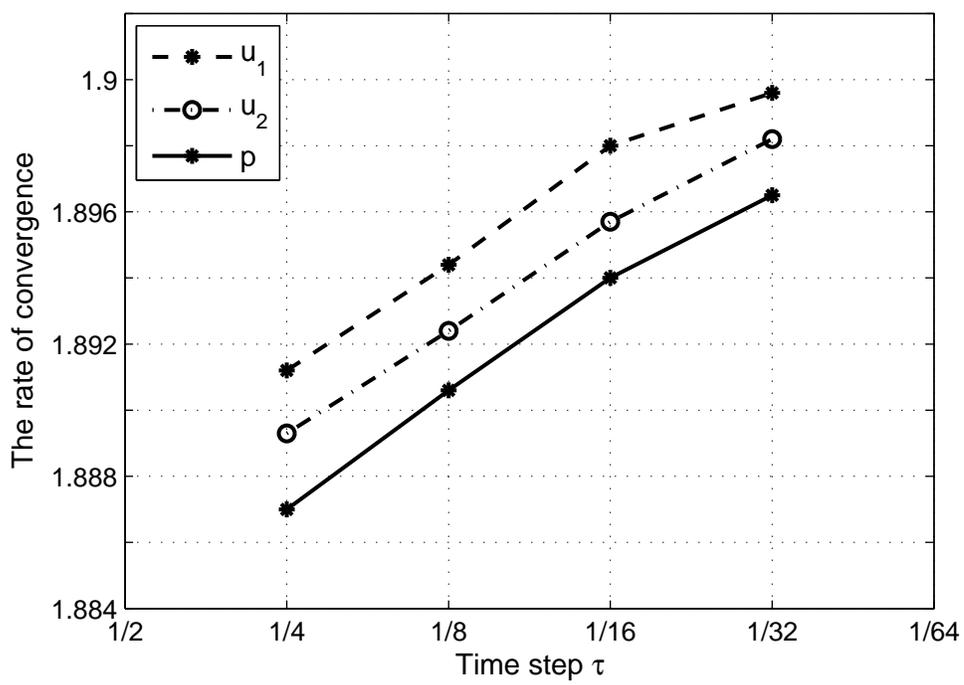
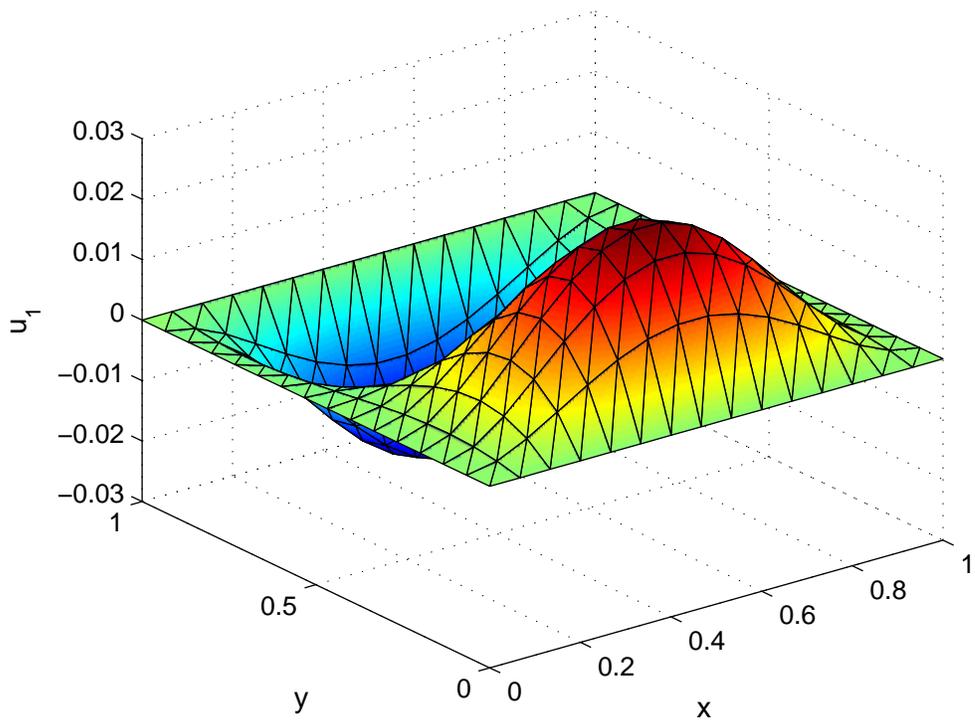
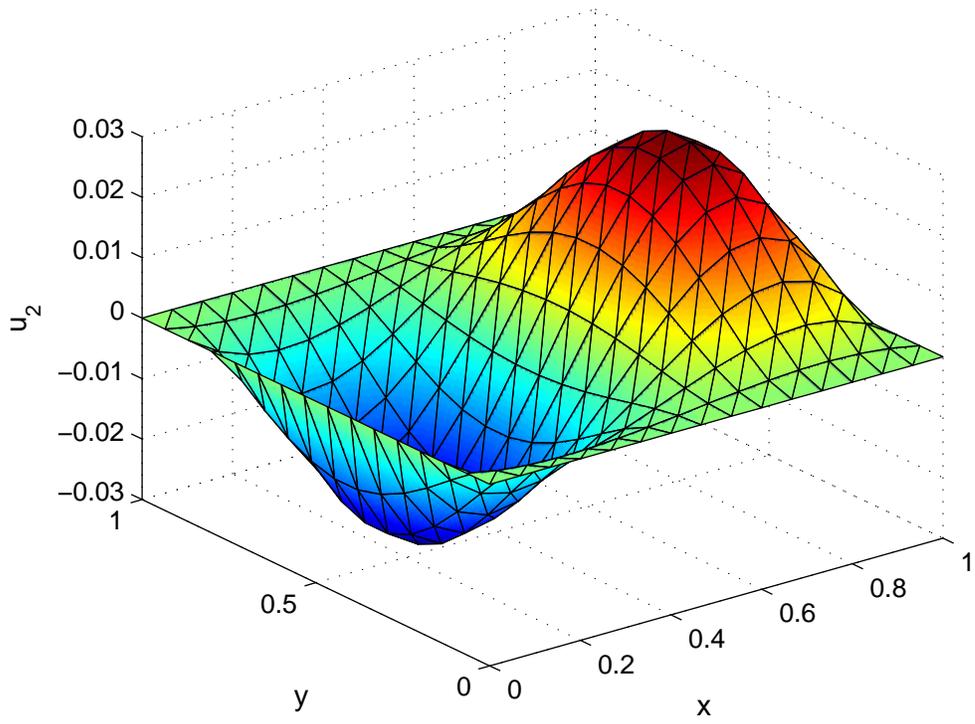
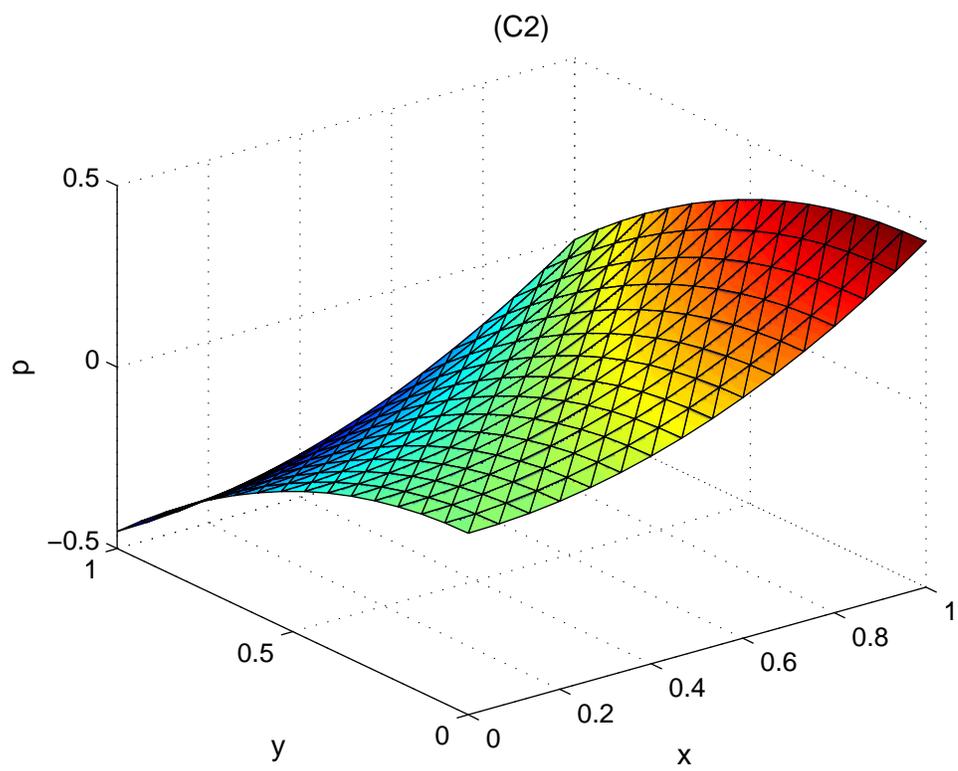
\end{document}